\documentclass[12pt,english]{article}
\usepackage[margin=0.75in]{geometry}
\usepackage[utf8]{inputenc}
\usepackage[english]{babel}
\usepackage{amsthm}
\usepackage{graphicx}
\usepackage{url}
\usepackage{hyperref}
\usepackage{amsmath}
\usepackage{listings}
\usepackage{amsfonts}
\usepackage{xcolor}
\usepackage{float}
\newtheorem{theorem}{Theorem}

\def\Xint#1{\mathchoice
   {\XXint\displaystyle\textstyle{#1}}%
   {\XXint\textstyle\scriptstyle{#1}}%
   {\XXint\scriptstyle\scriptscriptstyle{#1}}%
   {\XXint\scriptscriptstyle\scriptscriptstyle{#1}}%
   \!\int}
\def\XXint#1#2#3{{\setbox0=\hbox{$#1{#2#3}{\int}$}
     \vcenter{\hbox{$#2#3$}}\kern-.5\wd0}}
\def\dashint{\Xint-}

\begin{document}

\title{On the series expansion of the prime zeta function about $s=1$ and its coefficients}

\author{Artur Kawalec}

\date{}
\maketitle

\begin{abstract}
In this article, we derive a series expansion of the prime zeta function about the $s=1$ logarithmic singularity and prove general formula for its expansion coefficients, which is similar to the Stieltjes expansion coefficients for the Riemann zeta function. These results can also be viewed as a generalization of Mertens's Theorems to higher order. We also numerically verify and compute the presented formulas to high precision for several test cases.
\end{abstract}

\section{Introduction}
Let $p=\{2,3,5,7,\ldots\}$ be a sequence of primes and the prime zeta function

\begin{equation}\label{eq:1}
P(s)=\sum_{p}\frac{1}{p^s}
\end{equation}
is defined as a generalized zeta series over primes, which converges absolutely for $\Re(s)>1$. Through the M\"obius inversion formula as

\begin{equation}\label{eq:1}
P(s)=\sum_{n=1}^{\infty}\frac{\mu(n)}{n}\log\zeta(ns)
\end{equation}
the prime zeta can be further continued for $\Re(s)>0$. From this it appears that the prime zeta has logarithmic singularities (branch points) at $s=\frac{1}{n}$ and $\tfrac{\rho}{n}$ for all $n\geq 1$ whenever $\mu(n)\neq 0$ ($\rho$ is a non-trivial root of the zeta function). As a result, the prime zeta is a multi-valued function. And also, there is a natural boundary line at $\Re(s)=0$, after which $P(s)$ cannot be further continued to $\Re(s)\leq 0$ because of a dense accumulation of singularities that form near $\Re(s)=0$ [10, p. 215]. That's why branch points due to trivial zeros of the zeta function don't even enter $P(s)$.

In this article, we define a branch of $P(s)$ by the series expansion about the $s=1$ singularity as

\begin{equation}\label{eq:1}
P(s) =\log\left(\frac{1}{s-1}\right)+\sum_{n=0}^{\infty}\alpha_{n}\frac{(s-1)^n}{n!}
\end{equation}
The regular part of this series is analytic on a disk at $s=1$ with radius $R=\frac{1}{2}$, but because of the log term, we can define a branch cut on the line $(\tfrac{1}{2},1]$ inside the disk.

The main result of this article is we derive a general log-limit formula for these expansion coefficients as

\begin{theorem}
\label{thm:fta}
For $n=0$, the coefficient is:
\begin{equation}\label{eq:1}
\alpha_0 =\lim_{x\to\infty} \Bigg\{\sum_{p\leq x}\frac{1}{p}-\log\log(x)\Bigg\}-\gamma
\end{equation}
where $\gamma$ is Euler's constant, and for $n\geq 1$ the $n$th coefficient is:

\begin{equation}\label{eq:1}
\alpha_n =\lim_{x\to\infty} (-1)^n \Bigg\{\sum_{p\leq x}\frac{\log^n (p)}{p}-\frac{\log^n (x)}{n}\Bigg\}
\end{equation}

\end{theorem}

This limit formula is similar to the formula for Stieltjes constants in [1, p. 807], but it converges very slowly. But fortunately, the other representation by equation (2) involving the M\"obius function has much faster convergence, which we used to compute $\alpha_n$ to high precision in Table $1$.  So based on that these expansion coefficients can also be expressed by

\begin{theorem}
\label{thm:fta}
For $n\geq 0$ we have

\begin{equation}\label{eq:1}
\alpha_n=g_n+\frac{d^{n}}{ds^{n}}\Bigg[\sum_{k=2}^{\infty}\frac{\mu(k)}{k}\log\zeta(ks)\Bigg]_{s\to 1}
\end{equation}
\end{theorem}
\noindent where $g_n$ are the expansion coefficients of $\log\zeta(s)$ about $s=1$ (see Voros [11, p. 25]). The first few $g_n$ coefficients can be generated as

\begin{equation}\label{eq:1}
\begin{aligned}
g_0 &= 0, \\
g_1 &= \gamma, \\
g_2 &= -\gamma^2-2\gamma_1, \\
g_3 &= 2\gamma^3+6\gamma\gamma_1+3\gamma_2, \\
g_4 &= -6\gamma^4-12\gamma^2_1-24\gamma^2\gamma_1-12\gamma\gamma_2-4\gamma_3, \\
g_5 &= 120\gamma^3\gamma_1+120\gamma\gamma_1^2+60\gamma^2\gamma_2+60\gamma_1\gamma_2+20\gamma\gamma_3+5\gamma_4+24\gamma^5,\\
g_6 &= -720\gamma^4\gamma_1-1080\gamma^2\gamma_1^2-240\gamma_1^3-360\gamma^3\gamma_2-720\gamma\gamma_1\gamma_2-90\gamma_2^2+\\
&\frenchspacing-120\gamma^2\gamma_3-120\gamma_1\gamma_3-30\gamma\gamma_4-6\gamma_5-120\gamma^6, \\
\end{aligned}
\end{equation}
and so on, where they are expressed in terms of Stieltjes constants. These $g_n$ coefficients and their variations (often written by $\eta_n$) can also be found in Ma\'{s}lanka [6] and Bombieri-Lagarias [2] as being defined coefficients for $\zeta^{\prime}(s)/\zeta(s)$ function. They are generated by a recurrence relation of Coffey [3, p. 532] in terms of Stieltjes constants. The proof of Theorem 2 is to subtract the log singularity from $P(s)$ and differentiating as
\begin{equation}\label{eq:1}
\alpha_n=\frac{d^n}{ds^n}\Bigg[P(s) - \log\left(\frac{1}{s-1}\right)\Bigg]_{s\to 1}
\end{equation}
and substituting (2) for $P(s)$, since the log-zeta has similar series expansion
\begin{equation}\label{eq:1}
\log\zeta(s) =\log\left(\frac{1}{s-1}\right)+\sum_{n=0}^{\infty}g_{n}\frac{(s-1)^n}{n!}
\end{equation}
with $g_n$ coefficients. This log-zeta series is analytic on a disk centered at $s=1$ but with radius of convergence $R=3$ due to being limited by the first trivial zero of the zeta at $-2$, and so, because of the log term we can define a branch cut $(-2,1]$ for (9).

The first few cases of Theorem $1$ and $2$ are very well-known in the literature, i.e, the Mertens's Theorems [12]. For $n=0$ case is the harmonic sum of prime numbers

\begin{equation}\label{eq:1}
\begin{aligned}
M&=\lim_{x\to\infty} \Bigg\{\sum_{p\leq x}\frac{1}{p}-\log\log(x)\Bigg\}\\[1.2em]
&=0.26149721284764278375\ldots
\end{aligned}
\end{equation}
where the constant $M$ is known as the Meissel-Mertens's first constant (OEIS A077761), see also [7]. The relation

\begin{equation}\label{eq:1}
\alpha_0=M-\gamma
\end{equation}
shows that $\alpha_0$ is offset by Euler's constant

\begin{equation}\label{eq:1}
\gamma=0.5772156649015328606\ldots
\end{equation}
From (6) the value of
\begin{equation}\label{eq:1}
\begin{aligned}
\alpha_0&=\sum_{n=2}^{\infty}\frac{\mu(n)}{n}\log\zeta(n) \\[1.2em]
&= -0.31571845205389007685\ldots
\end{aligned}
\end{equation}
is the constant (OEIS A143524), which can be computed to high precision using this series involving the M\"obius function.

For $n=1$ the limit formula is

\begin{equation}\label{eq:1}
\begin{aligned}
\alpha_1 &=\lim_{x\to\infty} - \Bigg\{\sum_{p\leq x}\frac{\log(p)}{p}-\log (x)\Bigg\}  \\[1.2em]
&=  \gamma+\sum_{n=2}^{\infty}\mu(n)\frac{\zeta^{\prime}(n)}{\zeta(n)}\\[1.2em]
&= 1.33258227573322088176\ldots
\end{aligned}
\end{equation}
which is result of Mertens's Theorem (OEIS A083343) [12]. And for $n=2$ the formula is

\begin{equation}\label{eq:1}
\begin{aligned}
\alpha_2 &=\lim_{x\to\infty} \Bigg\{\sum_{p\leq x}\frac{\log^2(p)}{p}-\frac{\log^2(x)}{2}\Bigg\}  \\[1.2em]
&=   -\gamma^2-2\gamma_1+\sum_{n=2}^{\infty}\mu(n)n\left[-\frac{\zeta^{\prime}(n)^2}{\zeta(n)^2}+\frac{\zeta^{\prime\prime}(n)}{\zeta(n)}\right] \\[1.2em]
&= -2.55510761544644523959\ldots
\end{aligned}
\end{equation}
but such a constant doesn't seem to be computed in the literature. And one more example, for $n=3$ case we generate

\begin{equation}\label{eq:1}
\begin{aligned}
\alpha_3 &=\lim_{x\to\infty} -\Bigg\{\sum_{p\leq x}\frac{\log^3(p)}{p}-\frac{\log^3(x)}{3}\Bigg\}  \\[1.2em]
&= 2\gamma^3+6\gamma\gamma_1+3\gamma_2+\sum_{n=2}^{\infty}\mu(n)n^2\left[2\frac{\zeta^{\prime}(n)^3}{\zeta(n)^3}-3\frac{\zeta^{\prime}(n)\zeta^{\prime\prime}(n)}{\zeta(n)^2}+
\frac{\zeta^{\prime\prime\prime}(n)}{\zeta(n)}\right] \\[1.2em]
&= 10.25382709691100753877\ldots
\end{aligned}
\end{equation}
and so on. These $n$th log-zeta derivatives generate long such sequences of derivative terms in the series. In Table $1$ we computed the first $10$ coefficients numerically to high precision as a reference.

We also derived an alternate formula for these expansion coefficients during proof of Theorem $1$ as follows

\begin{theorem}
\label{thm:fta}
\begin{equation}\label{eq:1}
\alpha_n=(-1)^{n}\int_{1}^{\infty}\frac{\log^n(t)-n\log^{n-1}(t)}{t^2}[\pi(t)-\operatorname{li}(t)]\, dt
\end{equation}
\end{theorem}
\noindent
valid for $n\geq 0$ as shown in equation (34) and the logarithmic integral $\operatorname{li}(x)$ is defined in (21). 

We now turn our attention to proof of Theorem $1$ in next Section.

\begin{table}[hbt!]
\caption{Reference table for $\alpha_n$ coefficients}
\centering % used for centering table
\begin{tabular}{| c | c |} % centered columns (3 columns)
\hline
$\boldsymbol{n}$  & $\boldsymbol{\alpha_n}$ \textbf{(also in 30 digits)}\\ [0.5ex]
\hline
0 &-0.315718452053890076851085251473 \\
\hline
1 & 1.332582275733220881765828776071 \\
\hline
2 &-2.555107615446445239595583797989  \\
\hline
3 & 10.2538270969110075387787767411 \\
\hline
4 & -59.3323979717972728673195290222   \\
\hline
5 &  453.624590860932484915158069802  \\
\hline
6 & -4359.12496004203984785669925342  \\
\hline
7 &  50684.8409784215596972318317143   \\
\hline
8 & -692706.773919572383426686564824  \\
\hline
9 &  10884508.6063445498810870428549   \\
\hline
10& -193290090.992897724732297255085  \\
\hline
\end{tabular}
\label{table:nonlin} % is used to refer this table in the text
\end{table}

\section{Proof of Theorem 1}

To show that, let us recap some basic definitions. The prime counting function is defined by

\begin{equation}\label{eq:1}
\pi(x)=\sum_{p\leq x} 1
\end{equation}
for positive integer argument $x>0$, and by the average value

\begin{equation}\label{eq:1}
\pi(x)=\frac{1}{2}\Bigg[\sum_{p<x} 1 +\sum_{p\leq x} 1\Bigg]
\end{equation}
for positive real argument $x>0$. This means that its value is defined as an average of the two sides of the step whenever $x$ is integer. The asymptotic formula is

\begin{equation}\label{eq:1}
\pi(x)=\operatorname{li}(x)+f(x)
\end{equation}
where the logarithmic integral is defined by the P.V. integral function

\begin{equation}\label{eq:1}
\operatorname{li}(x)=\dashint_{0}^{x}\frac{1}{\log t} dt
\end{equation}
for real $x>0$. One also has the relation
\begin{equation}\label{eq:1}
\frac{d}{dx}\operatorname{li}(x)=\frac{1}{\log x}.
\end{equation}
And as for the remainder $f(x)$ the best estimate is
\begin{equation}\label{eq:1}
f(x)=O(\sqrt{x}\log x)
\end{equation}
assuming (RH). Also, the $m^{th}$ derivatives of (1) are given by
\begin{equation}\label{eq:1}
(-1)^mP^{(m)}(s)=\sum_{p}\frac{\log^m(p)}{p^s}
\end{equation}
for $\Re(s)>1$.

So we proceed proving equations (3-5) by carrying out appropriate Stieltjes integration. For the $m=0$ case (see Pain [9]) we write

\begin{align}
\sum_{p\leq x}\frac{1}{p} &= \int_{1+\delta}^{x} \frac{1}{t}d\pi(t)
= \int_{1+\delta}^{x} \frac{1}{t}d\operatorname{li}(t) + \int_{1+\delta}^{x} \frac{1}{t}df(t) \notag \\[1.2em]
&= \int_{1+\delta}^{x} \frac{1}{t\log(t)}\, dt
+ \left[ \frac{1}{t}f(t) \right]_{1+\delta}^{x} - \int_{1+\delta}^{x} \frac{d}{dt}\left(\frac{1}{t}\right)f(t)\, dt \notag \\[1.2em]
&= \Bigg[\log\log(x) -\log\log(1+\delta)\Bigg] + \Bigg[\frac{1}{x}f(x)-f(1+\delta)\Bigg]+\int_{1}^{x}\frac{1}{t^2}f(t)\, dt
\end{align}
Since the lower boundary term diverges at $t=1$, we take the limit $t=1+\delta$ as $\delta\to 0^{+}$. From (19) we have that
\begin{equation}\label{eq:1}
f(1+\delta)=-\operatorname{li}(1+\delta)
\end{equation}
and because $\pi(1^{-})=0$, $\pi(1)=0$ and $\pi(1^{+})=0$ we take the limit from positive side. The series expansion of log integral is

\begin{equation}\label{eq:1}
\operatorname{li}(x)=\gamma+\log\log(x)+\sum_{k=1}^{\infty}\frac{[\log(x)]^k}{k\:k!}
\end{equation}
so that near $1$ we have
\begin{equation}\label{eq:1}
\operatorname{li}(1+\delta)=\gamma+\log\log(1+\delta)+o(1)
\end{equation}
and as the divergent parts of the dominant and error terms cancel,  we get
\begin{equation}\label{eq:1}
\sum_{p\leq x}\frac{1}{p} = \log\log(x)+\gamma+\alpha_0+O\left(\frac{\log x}{\sqrt{x}}\right)
\end{equation}
And the integral converges to the $0$th order coefficient
\begin{equation}\label{eq:1}
\alpha_0=\int_{1}^{\infty}\frac{1}{t^2}f(t)dt.
\end{equation}

Now for the $m\geq 1$ case, we insert the $\log^m(t)$ factor in the integral and redo the calculation as follows:
\begin{align}
\sum_{p\leq x}\frac{\log^m(p)}{p} &= \int_{1}^{x} \frac{\log^m(t)}{t}d\pi(t)
= \int_{1}^{x} \frac{\log^m(t)}{t}d\operatorname{li}(t) + \int_{1}^{x} \frac{\log^m(t)}{t}df(t) \notag \\[1.2em]
&= \int_{1}^{x} \frac{\log^{m-1}(t)}{t}\, dt
+ \left[ \frac{\log^m(t)}{t}f(t) \right]_{1}^{x} - \int_{1}^{x} \frac{d}{dt}\left[\frac{\log^{m}(t)}{t}\right]f(t)\, dt \notag \\[1.2em]
&= A(x) + \frac{\log^m(x)}{x}f(x)-\int_{1}^{x}\frac{m\log^{m-1}(t)-\log^m(t)}{t^2}f(t)\, dt
\end{align}
and this integral identity can be generated as

\begin{equation}\label{eq:1}
\begin{aligned}
A(x)&=\int_{1}^{x} \frac{\log^{m-1}(t)}{t}\, dt\\[1.2em]
&=\frac{\log^{m}(x)}{m}
\end{aligned}
\end{equation}
for $m\geq 1$. Then, it is readily seen that the limit

\begin{equation}\label{eq:1}
\sum_{p\leq x}\frac{\log^m(p)}{p}-\frac{\log^{m}(x)}{m}=(-1)^m \alpha_m+O\left(\frac{\log^{m+1}(x)}{\sqrt{x}}\right)
\end{equation}
exists as $x\to\infty$ and converges to a constant

\begin{equation}\label{eq:1}
\alpha_m=(-1)^{m}\int_{1}^{\infty}\frac{\log^m(t)-m\log^{m-1}(t)}{t^2}f(t)\, dt
\end{equation}
for $m\geq 1$ by the estimate (23). Now the next step is to connect this constant to the series expansion in (3) as follows:

\begin{equation}\label{eq:1}
\begin{aligned}
P(s)
&= \int_{1+\delta}^{\infty}\frac{1}{t^{s}}\, d\pi(t) = \int_{1+\delta}^{\infty} \frac{1}{t^{s}} d\operatorname{li}(t)\,  + \int_{1+\delta}^{\infty} \frac{1}{t^{s}} df(t) \\[1.2em]
&=\int_{1+\delta}^{\infty} \frac{1}{t^{s}\log t} dt\,  + \Biggl[\frac{f(t)}{t^s}\Biggr]_{1+\delta}^\infty
 + s\int_{1+\delta}^\infty t^{-s-1}f(t)\, dt \\[1.2em]
\end{aligned}
\end{equation}
for $\Re(s)>1$ as $\delta\to 0^{+}$, where here again, the integrand diverges at $t=1$ where it has a simple pole. The anti-derivative is related to the logarithmic or exponential integral function (21), but we express the anti-derivative in terms of the exponential integral as

\begin{equation}\label{eq:1}
\int_{1+\delta}^{\infty} \frac{1}{t^{s}\log t} dt=\operatorname{E_1}[(s-1)\log(1+\delta)]
\end{equation}
by substitution $u=\log t$ and $t=e^u$ and $dt=e^u du$ yields the definition 

\begin{equation}\label{eq:1}
\operatorname{E_1}(z)=\int_{z}^{\infty}\frac{e^{-t}}{t}dt, \quad |\arg(z)|<\pi
\end{equation}
of the exponential integral valid for all complex $z$, but with a branch cut on the negative real axis $(-\infty,0]$ as defined in [1, p. 228]. The function $\operatorname{E_1}(z)$ also admits the series expansion

\begin{equation}\label{eq:1}
\operatorname{E_1}(z)=-\gamma-\log(z)-\sum_{k=1}^{\infty}\frac{(-z)^k}{k\:k!}
\end{equation}
And so, the result of (36) has the expansion
\begin{equation}\label{eq:1}
\begin{aligned}
\operatorname{E_1}[(s-1)\log(1+\delta)]&=-\gamma-\log [(s-1)\log(1+\delta)]-\sum_{k=1}^{\infty}\frac{[-(s-1)\log(1+\delta)]^k}{k\:k!} \\[1.2em]
&=-\gamma-\log(s-1)-\log\log(1+\delta)+o(1)
\end{aligned}
\end{equation}
where it diverges as $\delta\to 0^{+}$. The logarithm of product split is justified because $\arg(1+\delta)\to 0$.  Then the second term arising from the boundary condition of the second integral in (35) also diverges

\begin{equation}\label{eq:1}
\begin{aligned}
f(1+\delta)&=-\operatorname{li}(1+\delta)\\[1.2em]
&=-\gamma-\log\log(1+\delta)+o(1)
\end{aligned}
\end{equation}
As a result, on combining (39) and (40) one takes the limit

\begin{equation}\label{eq:1}
\lim_{\delta\to 0^{+}}\Bigg[\operatorname{E_1}((s-1)\log(1+\delta))-f(1+\delta)\Bigg]=\log\left(\frac{1}{s-1}\right)
\end{equation}
there is a cancelation of two divergent quantities for $\Re(s)>1$, where we recover the log term, but then the log term analytically extends to larger domain except for the branch cut $(-\infty,1]$. And so, this results in
\begin{equation}\label{eq:1}
P(s) = \log\left(\frac{1}{s-1}\right) + s\int_1^\infty t^{-s-1}f(t)\, dt
\end{equation}
the integral in (42) converges for $\Re(s)>\tfrac{1}{2}$ when one inserts (23) assuming (RH). As a result, it defines an analytical continuation of $P(s)$ for $\Re(s)>\tfrac{1}{2}$ with a branch cut $(\frac{1}{2},1]$.

It now remains to show that this integral admits the series expansion about $s=1$, so we consider

\begin{equation}\label{eq:1}
\int_{1}^{\infty} t^{-s-1} f(t)\, dt
\end{equation}
which converges for $\Re(s)>\tfrac{1}{2}$ assuming (RH). One then generates the exp-log expansion as
\begin{equation}\label{eq:1}
\int_{1}^{\infty} t^{-2} e^{-(s-1) \log t}f(t) \, dt = \sum_{j=0}^{\infty} c_j \frac{(s-1)^j}{j!}
\end{equation}
where its expansion coefficients can be read off as
\begin{equation}\label{eq:1}
c_j = (-1)^j\int_{1}^{\infty} \frac{\log^j(t)}{t^{2}}f(t)  \, dt
\end{equation}
but such a series is only limited to $|s-1|<1$ domain. However, on writing integral term in (42) in $s-1$ domain, then one needs to consider this form instead

\begin{equation}\label{eq:1}
\begin{aligned}
s\int_{1}^{\infty} t^{-s-1} f(t)\, dt &= (s-1)\int_{1}^{\infty} t^{-s-1} f(t)\, dt+\int_{1}^{\infty} t^{-s-1} f(t)\, dt \\[1.2em]
&= \int_{1}^{\infty} [(s-1)t^{-s-1} f(t)+t^{-s-1} f(t)]\, dt \\[1.2em]
\end{aligned}
\end{equation}
where we can generate a new series expansion as
\begin{equation}\label{eq:1}
\begin{aligned}
s\int_{1}^{\infty} t^{-s-1} f(t)\, dt&=\sum_{j=0}^{\infty}(j+1)c_j\frac{(s-1)^{j+1}}{(j+1)!}+\sum_{j=0}^{\infty}c_j\frac{(s-1)^{j}}{j!}\\[1.2em]
&=\sum_{m=1}^{\infty}mc_{m-1}\frac{(s-1)^{m}}{m!}+c_0+\sum_{m=1}^{\infty}c_m\frac{(s-1)^{m}}{m!}\\[1.2em]
&=c_0+\sum_{m=1}^{\infty}\left[mc_{m-1}+c_m\right]\frac{(s-1)^{m}}{m!}\\[1.2em]
&=\sum_{j=0}^{\infty} y_j \frac{(s-1)^j}{j!} \\[1.2em]
\end{aligned}
\end{equation}
by combining the expansions in (46). Here we've shifted index variable $m=j+1$ in first series in (47) and re-labeled $m=j$ in second series in (47). The new expansion coefficients are then

\begin{equation}\label{eq:1}
y_m=\begin{cases}
c_0 & \text{for } m = 0 \\
mc_{m-1}+c_m & \text{for } m \geq 1
\end{cases}
\end{equation}
where they are expressed by $c_n$ coefficients, and by (45) these coefficients can be expressed by the integral

\begin{equation}\label{eq:1}
y_m=(-1)^m\int_{1}^{\infty}\frac{\log^m(t)-m\log^{m-1}(t)}{t^2}f(t)\, dt
\end{equation}
which matches the result in (34), hence for the $m\geq 0$ we have $\alpha_m=y_m$. Hence this completes the proof.

\texttt{Email: art.kawalec@gmail.com}

\end{document}